\documentclass{amsart}%
\usepackage{amsfonts}
\usepackage{amsmath}
\usepackage{amssymb}
\usepackage{graphicx}%
\newtheorem{theorem}{Theorem}[section]

\theoremstyle{definition}

\theoremstyle{remark}
\newtheorem{remark}[theorem]{Remark}
\numberwithin{equation}{section}
\theoremstyle{plain}
\newtheorem{acknowledgement}{Acknowledgement}

\copyrightinfo{2015}{David Glickenstein and Liang Wu}
\begin{document}
\title[Solitons for RG-2 flow]{Soliton metrics for two-loop renormalization group flow on 3D unimodular Lie groups}
\author{David Glickenstein}
\address{University of Arizona\\
Tucson, AZ 85721}
\email{glickenstein@math.arizona.edu}
\author{Liang Wu}

\begin{abstract}
The two-loop renormalization group flow is studied via the induced bracket
flow on 3D unimodular Lie groups. A number of steady solitons are found. Some
of these steady solitons come from maximally symmetric metrics that are
steady, shrinking, or expanding solitons under Ricci flow, while others are
not obviously related to Ricci flow solitons.

\end{abstract}
\maketitle

\section{Introduction}

In this note we look into perturbation of the Ricci flow on three-dimensional
unimodular Lie groups. The Ricci flow on three-dimensional Lie groups has been
well-studied, both in forward and backward time \cite{IJ}\cite{KM}%
\cite{CS}\cite{GP}\cite{CGS}, as well as other geometric flows such as cross
curvature flow \cite{CNS}\cite{CS2} and other flows, e.g., \cite{KST}%
\cite{BBLP}\cite{PPS}. We will use the notations $E\left(  2\right)  $ for the
universal covering of the group of isometries of the Euclidean plane,
$H\left(  3\right)  $ for the three-dimensional Heisenberg group, $E\left(
1,1\right)  $ for the group of isometries of the Minkowski plane,
$\operatorname{SU}\left(  2\right)  $ for the special unitary group, and
$\operatorname{SL}\left(  2,\mathbb{R}\right)  $ for the universal covering of
the special linear group (this is nonstandard notation but simplifies the
text). Many of these groups admit maximally symmetric Thurston geometries that
are special left-invariant metrics (actually, one- and two-parameter families
of metrics) with additional symmetries:\ the flat metric on $E\left(
2\right)  ,$ $\operatorname{Nil}$ on $H\left(  3\right)  ,$
$\operatorname{Sol}$ on $E\left(  1,1\right)  ,$ the round sphere
$\mathbb{S}^{3}$ metrics on $\operatorname{SU}\left(  2\right)  ,$ and the
maximally symmetric metrics on $\operatorname{SL}\left(  2,\mathbb{R}\right)
.$ It is well-known that many of these manifolds admit steady Ricci solitons
(the flat metrics), shrinking Ricci solitons (the round sphere metrics), and
expanding Ricci solitons ($\operatorname{Nil}$, $\operatorname{Sol}$). It is
notable, however, that not all of Thurston's geometries are represented among
these generalized fixed points; in particular, there is no soliton on
$\operatorname{SL}\left(  2,\mathbb{R}\right)  .$ Also note that the case of
expanding solitons is a bit misleading, as the solitons do not exist on
compact quotients, but instead lead to collapsing of closed manifolds with
these geometries.

One may ask if we can deform the Ricci flow to another flow that has
additional fixed points, steady solitons, or other special solutions. One
candidate for such a perturbation is the two-loop renormalization group flow
(RG-2 Flow), introduced by physicists and studied more recently in works such
as \cite{GGI2}\cite{CM}\cite{GGI3}. The flow is compactly written as
\begin{equation}
\frac{\partial}{\partial t}g=-2\operatorname{Rc}-\frac{\alpha}{2}%
\operatorname{Rm}^{2} \label{eq: RG2}%
\end{equation}
where $\alpha$ is a constant.

In this note, we describe several steady soliton solutions to RG-2 flow
depending on the value of the parameter $\alpha$. It was recently brought to
the authors' attention that a similar classification of solitons for RG-2 flow
was found concurrently by T. Wears in \cite{We} by finding derivations on the
Lie algebras. Our approach is instead to look at the flow of Lie bracket
coefficients on an orthnormal frame, which has the advantage of making it
easier (and quicker) to find soliton metrics but the disadvantage of not
having an explicit family of automorphisms.

\begin{acknowledgement}
The first author would like to thank Tracy Payne and Christine Guenther, whose
research and informal discussions paved the way for this work. The authors
were partially funded by NSF grant DMS 0748283.
\end{acknowledgement}

\section{RG-2 bracket flow}

Instead of looking at the RG-2 flow on the metric coefficients, we prefer to
work on a left-invariant orthonormal frame. It is well known \cite{Mil} that
three-dimensional unimodular Lie algebras can be written in terms of an
orthonormal frame with at most three nonzero structure constants:%
\begin{equation}
\lbrack e_{2},e_{3}]=a_{1}e_{1,}\quad\lbrack e_{3},e_{1}]=a_{2}e_{2}%
,\quad\text{and }[e_{1},e_{2}]=a_{3}e_{3}, \label{eq:1}%
\end{equation}
for some constants $a_{1},a_{2},a_{3}\in\mathbb{R}.$ The structure constants
determine the left-invariant metric (see e.g., \cite{LL}), and so any flow of
metrics induces a flow of structure constants. This induced flow is called the
corresponding bracket flow. The bracket flow for Ricci flow has been studied
in a number of contexts. Note that the sign of the structure constants
$a_{1},a_{2},a_{3}$ determine which unimodular Lie group the algebra
corresponds to among $\mathbb{R}^{3},$ $H\left(  3\right)  ,$ $E\left(
3\right)  ,$ $E\left(  1,1\right)  ,$ $\operatorname{SL}\left(  2,\mathbb{R}%
\right)  ,$ $\operatorname{SU}\left(  2\right)  .$

A calculation along the lines of that in \cite{GP} leads to the following
system on the structure constants corresponding to the RG-2 flow:%
\begin{align*}
\dfrac{da_{1}}{dt}  &  =\frac{a_{1}}{2}\left[  (a_{2}-a_{3})^{2}-3a_{1}%
^{2}+2a_{1}a_{2}+2a_{1}a_{3}\right]  \left\{  1+\frac{\alpha}{8}\left[
(a_{2}-a_{3})^{2}-3a_{1}^{2}+2a_{1}a_{2}+2a_{1}a_{3}\right]  \right\} \\
\dfrac{da_{2}}{dt}  &  =\frac{a_{2}}{2}\left[  (a_{1}-a_{3})^{2}-3a_{2}%
^{2}+2a_{1}a_{2}+2a_{2}a_{3}\right]  \left\{  1+\frac{\alpha}{8}\left[
(a_{1}-a_{3})^{2}-3a_{2}^{2}+2a_{1}a_{2}+2a_{2}a_{3}\right]  \right\} \\
\dfrac{da_{3}}{dt}  &  =\frac{a_{3}}{2}\left[  (a_{1}-a_{2})^{2}-3a_{3}%
^{2}+2a_{1}a_{3}+2a_{2}a_{3}\right]  \left\{  1+\frac{\alpha}{8}\left[
(a_{1}-a_{2})^{2}-3a_{3}^{2}+2a_{1}a_{3}+2a_{2}a_{3}\right]  \right\}  .
\end{align*}
Note that when $\alpha=0,$ the bracket formulation of the Ricci flow is recovered.

An important property of the bracket flow is that steady algebraic solitons
are precisely the metrics that are fixed under the bracket flow (see, e.g.,
\cite{LL}). A steady algebraic soliton is a left invariant metric $g\left(
t\right)  $ such that
\[
g\left(  t\right)  =\phi_{t}^{\ast}g_{0}%
\]
where $\phi_{t}$ is a one-parameter family of automorphisms of the group (the
use of \textquotedblleft algebraic\textquotedblright\ refers to the assumption
that $\phi_{t}$ are automorphisms, not just diffeomorphisms). RG-2 flow
(\ref{eq: RG2}) is not scale invariant in the same way that Ricci flow is
scale invariant due to the fact that the constant $\alpha$ is not
dimensionless; usually we specify a value of $\alpha$ and that fixes a scale.
For this reason, we are mainly concerned with steady solitons.

\section{Solitons for RG-2 flow}

The main result is the following.

\begin{theorem}
The following are all steady solitons for the RG-2 flow, characterized by
their structure constants $\left(  a_{1},a_{2},a_{3}\right)  ,$ the sign of
$\alpha,$ their sectional curvatures (given as the three eigenvalues of the
Einstein tensors), and their Ricci curvatures. The maximally symmetric metrics
are described by their geometries.

\hspace*{-0.8in}%
\begin{tabular}
[c]{||l||l||l||l||l||l||}\hline
\textbf{Group} & $\left(  a_{1},a_{2},a_{3}\right)  $ & $\mathbf{\alpha}$ &
\begin{tabular}
[c]{l}%
\textbf{Sectional}\\
\textbf{curvatures}%
\end{tabular}
\textbf{ } &
\begin{tabular}
[c]{l}%
\textbf{Ricci}\\
\textbf{curvatures}%
\end{tabular}
& \textbf{Geometry}\\\hline\hline
$E\left(  2\right)  $ & $\left(  a,a,0\right)  ,$ $a\in\mathbb{R}$ & any &
$\left(  0,0,0\right)  $ & $\left(  0,0,0\right)  $ & flat\\\hline
$H\left(  3\right)  $ & $\left(  \pm\sqrt{\frac{3}{8\alpha}},0,0\right)  $ &
$\alpha>0$ & $\left(  \frac{3}{32\alpha},\frac{3}{32\alpha},-\frac{9}%
{32\alpha}\right)  $ & $\left(  \frac{3}{16\alpha},-\frac{3}{16\alpha}%
,-\frac{3}{16\alpha}\right)  $ & $\operatorname{Nil}$\\\hline
$E\left(  1,1\right)  $ & $\left(  \pm\sqrt{\frac{2}{-\alpha}},0,\mp
\sqrt{\frac{2}{-\alpha}}\right)  $ & $\alpha<0$ & $\left(  \frac{2}{\alpha
},\frac{2}{\alpha},-\frac{2}{\alpha}\right)  $ & $\left(  0,0,\frac{4}{\alpha
}\right)  $ & $\operatorname{Sol}$\\\hline
$\operatorname{SU}\left(  2\right)  $ & $\left(  \pm\sqrt{\frac{8}{\alpha}%
},\pm\sqrt{\frac{8}{\alpha}},\pm\sqrt{\frac{8}{\alpha}}\right)  $ & $\alpha>0$
& $\left(  \frac{2}{\alpha},\frac{2}{\alpha},\frac{2}{\alpha}\right)  $ &
$\left(  \frac{4}{\alpha},\frac{4}{\alpha},\frac{4}{\alpha}\right)  $ &
$\mathbb{S}^{3}\left(  \sqrt{\frac{\alpha}{2}}\right)  $\\\hline
$E\left(  1,1\right)  $ & $\left(  -\frac{3}{2\sqrt{\alpha}},0,\frac{1}%
{2\sqrt{\alpha}}\right)  $ & $\alpha>0$ & $\left(  \frac{1}{\alpha}%
,0,\frac{-2}{\alpha}\right)  $ & $\left(  \frac{1}{\alpha},-\frac{2}{\alpha
},-\frac{1}{\alpha}\right)  $ & \\\hline
$\operatorname{SU}\left(  2\right)  $ & $\left(  \frac{3}{\sqrt{-2\alpha}%
},\frac{3}{\sqrt{-2\alpha}},\frac{4}{\sqrt{-2\alpha}}\right)  $ & $\alpha<0$ &
$\left(  \frac{2}{\alpha},\frac{2}{\alpha},-\frac{6}{\alpha}\right)  $ &
$\left(  -\frac{4}{\alpha},-\frac{4}{\alpha},\frac{4}{\alpha}\right)  $ &
\\\hline
$\operatorname{SU}\left(  2\right)  $ & $\left(  \frac{3}{\sqrt{-2\alpha}%
},\frac{3}{\sqrt{-2\alpha}},\frac{2}{\sqrt{-2\alpha}}\right)  $ & $\alpha<0$ &
$\left(  -\frac{1}{2\alpha},-\frac{1}{2\alpha},\frac{3}{8\alpha}\right)  $ &
$\left(  -\frac{4}{\alpha},-\frac{4}{\alpha},-\frac{2}{\alpha}\right)  $ &
\\\hline
\end{tabular}

\end{theorem}

\begin{remark}
It is important to note that the first four metrics are solitons for Ricci
flow, correspond to a metric in the one-parameter family of maximally
symmetric metrics (which are Thurston geometries). Because they are maximally
symmetric, any metric starting in the one-parameter family remains in it.
These particular choices of initial metric for certain $\alpha$ are actually
fixed under RG-2, whereas they may be expanding or shrinking under Ricci flow.
In addition to the fixed points listed here, it is easy to use the bracket
flow to describe the dynamics under the system (\ref{eq:1}) of any maximally
symmetric metric of these forms even if they are not fixed; they correspond to
self-similar expanding or shrinking solutions.
\end{remark}

\begin{remark}
According to \cite{CM} and \cite{GGI3}, the RG-2 flow is parabolic if
$1+\alpha K>0$ for all sectional curvatures $K.$ The only solitons in the
table that satisfy this requirement are the first two, the fourth, and the last.
\end{remark}

The last three metrics are new compared to the work in \cite{GGI}, especially
the metric on $E\left(  1,1\right)  ,$ which does not satisfy the symmetry
assumption in that paper; quite possibly \cite{GGI} were aware of the last two
but did not consider these solutions since they require $\alpha<0,$ which is
not physical. As we are only considering the system as a perturbation of Ricci
flow, we have no problem considering $\alpha<0.$

\section{Sketch of proof and discussion}

The main idea is to start with the ansatz that there is a solution of one of
the following forms: $\left(  a\left(  t\right)  ,0,c\left(  t\right)
a\left(  t\right)  \right)  $ or $\left(  a\left(  t\right)  ,a\left(
t\right)  ,c\left(  t\right)  a\left(  t\right)  \right)  $ and derive the
equations to maintain these forms. One can then derive the equations for $c$
to be constant. Any such solution will be self-similar, and this can be used
to analyze self-similar solutions similar to those in the first four cases,
where we get shrinkers or expanders; in particular, we can find when such
solutions are fixed, which depends on the initial condition $a\left(
0\right)  .$ For the last three solutions, we cannot derive a solution where
$c$ is constant without $a$ being constant as well. The solutions where $a$
and $c$ are both constant correspond to algebraic steady soliton metrics.

This method was motivated by studying a related problem. In \cite{GP}, the
Ricci flow was considered up to scaling by looking at the ratio of structure
constants $m_{2}=\frac{a_{2}}{a_{1}}$ and $m_{3}=\frac{a_{3}}{a_{1}}.$ This
allowed a complete phase portrait of the system described up to scaling, and
hence in the $\left(  m_{2},m_{3}\right)  $-plane, expanding and shrinking
solitons are also fixed points. As previously mentioned, RG-2 flow does not
allow this analysis since it is not scale invariant, but one could still make
an attempt to consider the flow in the $\left(  m_{2},m_{3}\right)  $-plane,
getting the equations%
\begin{equation}%
\begin{array}
[c]{c}%
\frac{dm_{2}}{dt}=m_{2}\left(  1-m_{2}\right)  \left(  1+m_{2}-m_{3}\right)
\left[  1-\beta\left(  1+m_{3}-m_{2}\right)  \left(  1-m_{2}-m_{3}\right)
\right]  ,\\
\frac{dm_{3}}{dt}=m_{3}\left(  1-m_{3}\right)  \left(  1+m_{3}-m_{2}\right)
\left[  1-\beta\left(  1+m_{2}-m_{3}\right)  \left(  1-m_{2}-m_{3}\right)
\right]  ,
\end{array}
\label{eq:7}%
\end{equation}
where $\beta=\frac{\alpha a_{1}^{2}}{4}.$ Note that, in general, $\beta$ is
not constant since $\alpha$ is constant and $a_{1}$ is not. However, for
steady solitons of RG-2 flow, $\beta$ is constant and so this system gives
corresponding fixed points.

The system (\ref{eq:7}) can be solved completely in the sense that we can find
all the fixed points and analyze the phase plane analogously to \cite{GP}. The
fixed points are $\left(  0,0\right)  ,$ $\left(  0,\pm1\right)  ,$ $\left(
\pm1,0\right)  ,$ $\left(  1,1\right)  ,$ $\left(  0,1\pm\sqrt{\dfrac{1}%
{\beta}}\right)  ,$ $\left(  1\pm\sqrt{\dfrac{1}{\beta}},0\right)  ,$ $\left(
1,1\pm\sqrt{1+\dfrac{1}{\beta}}\right)  ,$ $\left(  1\pm\sqrt{1+\dfrac
{1}{\beta}},1\right)  ,\text{ and $\dfrac{1}{2}\left(  1-\dfrac{1}{\beta
},1-\dfrac{1}{\beta}\right)  $}.$ Notice that this system has fixed points
(depending on the value of $\beta$) in a number of regimes including
$\operatorname{SL}\left(  2,\mathbb{R}\right)  ,$ which does not have a fixed
point for the bracket flow of either the Ricci flow or RG-2 flow. However, it
is not clear how this flow corresponds to any geometric flow, since it would
require a time dependent $\alpha$, determined in a way that is not calculated
from a curvature.

\end{document}